\documentclass{amsart}
\usepackage{amssymb}
\usepackage{amsbsy}
\usepackage{upgreek}
\usepackage{hyperref}

\input xy
\xyoption{all}

\begin{document}

\renewcommand{\sc}{\scshape}
\newcommand{\D}{\mathcal{D}}
\newcommand{\G}{\mathbf{G}}
\renewcommand{\H}{\mathbf{H}}
\newcommand{\K}{\mathbf{K}}
\newcommand{\X}{\mathbf{X}}
\newcommand{\Y}{\mathbf{Y}}
\newcommand{\Z}{\mathbf{Z}}
\newcommand{\XI}{X_\infty}
\newcommand{\YI}{Y_\infty}
\newcommand{\F}{\mathbf{F}}
\newcommand{\f}{\mathbf{f}}
\newcommand{\g}{\mathbf{g}}
\newcommand{\p}{\mathbf{p}}
\newcommand{\q}{\mathbf{q}}
\newcommand{\Bpi}{\boldsymbol{\uppi}}
\newcommand{\Bphi}{\boldsymbol{\upvarphi}}
\newcommand{\Bpsi}{\boldsymbol{\uppsi}}
\newcommand{\Aut}{\mathrm{Aut}}
\newcommand{\VAut}{\mathrm{VAut}}
\newcommand{\newword}[1]{\textbf{#1}}

\title{Compact aspherical solenoids}

\author{James Belk}
\address{University of Glasgow, Glasgow, Scotland}
\email{jim.belk@glasgow.ac.uk}
\thanks{The first author was partially supported by an NSF Postdoctoral Research Fellowship.}

\author{Bradley Forrest}
\address{Stockton University \\ 101 Vera King Farris Drive,
Galloway, NJ 08205}
\email{bradley.forrest@stockton.edu}



\newtheorem{theorem}{Theorem}[section]
\newtheorem{lemma}[theorem]{Lemma}
\newtheorem{corollary}[theorem]{Corollary}
\newtheorem{proposition}[theorem]{Proposition}

\theoremstyle{definition}
\newtheorem{definition}[theorem]{Definition}
\newtheorem{example}[theorem]{Example}

\theoremstyle{remark}
\newtheorem{remark}[theorem]{Remark}
\newtheorem{assumptions}[theorem]{Assumptions}
\numberwithin{equation}{section}

\begin{abstract}
We consider compact, aspherical solenoids obtained as the inverse limit of a system of CW~complexes and covering maps.  This includes $P$-adic solenoids, as well as the universal hyperbolic solenoid of Teichm\"{u}ller theory.  Using ideas from shape theory, we classify maps between such solenoids up to homotopy, and we prove a Dehn-Nielsen-type theorem for self-homotopy equivalences of such a solenoid.  This generalizes a result of Odden regarding the universal hyperbolic solenoid.
\end{abstract}

\maketitle

\section{Introduction}

Classically, a \textit{solenoid} is the inverse limit of a sequence of circles
\[
S^1 \longleftarrow S^1  \longleftarrow  S^1  \longleftarrow  \cdots,
\]
where each map is a covering map \cite{vandantzig}.  Such spaces, now referred to as $P$-adic solenoids, can be described as intersections of solid tori in $\mathbb{R}^3$, and arise as strange attractors of certain dynamical systems.  The topology and dynamics of these solenoids have been studied by several authors \cite{charatonik-covarrubias} \cite{gumerov} \cite{takens}.

In 1965, McCord \cite{mccord} and Schori \cite{schori} generalized the notion of a solenoid to include any inverse limit of a sequence of ``nice'' spaces and regular covering maps.  Since then, these more general solenoids have been studied extensively \cite{clark-fokkink2} \cite{clark-fokkink} \cite{fokkink-oversteegen}.  In addition to their importance as strange attractors, solenoids of this type have played an important role in the topological theory of continua \cite{eda-mandic-matijevic} \cite{krupski-rogers} \cite{kwapisz}.  Such solenoids also play an increasing role in Teichm\"{u}ller theory, where the the universal hyperbolic solenoid defined by Sullivan \cite{sullivan} has become an important object of study \cite{biswas-nag-sullivan} \cite{markovic-saric}.

The motivation for this note comes from the work of Odden \cite{odden} on the mapping class group of the universal hyperbolic solenoid.  From our point of view, Odden's argument uses hyperbolic geometry in an essential way to prove what ought to be a theorem of topology.  In Section~8, we present a purely topological version of Odden's theorem that holds true for compact aspherical solenoids.  For solenoids in this class, the homotopy type is entirely determined by the fundamental pro-group, and homotopy classes of maps are in one-to-one correspondence with pro-group morphisms.

Our results depend heavily on certain theorems from shape theory, particularly the inverse system approach \cite{MaSe}.  Though we are not experts on this subject, we have found that many of the standard definitions and theorems of shape theory become simpler in the context of inverse systems of covers.  We have therefore included a basic exposition of the shape theory that we require, specialized to this context.  We hope this will be useful for researchers in dynamical systems, Teichm\"{u}ller theory, and other fields who are interested in solenoids but may be unfamiliar with shape theory.

This paper is organized as follows.  Section~2 introduces covering systems and solenoids, and reviews some basic examples.  Section~3 discusses morphisms between covering systems, and the resulting limit maps between solenoids.  In Section~4 we recall some theorems from shape theory on expansions, and apply these to our context. Section~5 introduces pro-groups and filtered groups, which are the algebraic counterparts to covering systems in the theory, and Section~6 introduces the fundamental pro-group of a solenoid.  In Section~7 we apply our theory to classify homotopy classes of maps between certain aspherical solenoids, and in Section~8 we apply these results to virtual automorphism groups, generalizing Odden's theorem.

Many of the results in Sections 2 through 6 are special cases of standard definitions and theorems from shape theory, and we make no claim of originality for anything but the context and style of exposition.

Most of the material in Sections~7~and~8 is new, including Theorems~\ref{thm:nielsen}, \ref{cor:homotopyclass}, \ref{cor:mappingclassgroup}, and \ref{thm:universalsolenoid} as well as most of the applications of these theorems discussed in the subsequent examples.

\subsection*{Acknowledgements}
The first draft of this manuscript appeared in 2010, but was never published.  The ideas here have since proven themselves useful in the study of mapping class groups of solenoids (see~\cite{bering-studenmund}).  We would like to thank Steven Hurder for encouraging us to prepare this revised version.

\section{Covering systems and solenoids}

We will use the term \newword{covering system} to refer to any inverse system of connected CW complexes and covering maps, indexed over a directed set.  For example, any inverse sequence of connected covers
\[
X_0 \;\overset{p_1}{\longleftarrow}\; X_1
\;\overset{p_2}{\longleftarrow}\; X_2
\;\overset{p_3}{\longleftarrow}\; \cdots
\]
of a CW complex $X_0$ is a covering system, indexed over the natural numbers.  In general, a covering system $\X$ consists of a collection of CW complexes $\{X_\alpha\}_{\alpha\in\D}$ indexed over a directed set~$\D$, with covering maps $p_{\alpha\beta}\colon X_\beta \to X_\alpha$ for all $\alpha, \beta \in \D$ with $\alpha\leq\beta$. These maps are required to commute in the sense that $p_{\alpha\beta} \circ p_{\beta\gamma} = p_{\alpha\gamma}$ for all $\alpha,\beta,\gamma\in\D$ with $\alpha \leq \beta \leq \gamma$.

\begin{remark}
We will assume throughout that all topological spaces come with a basepoint, and that all maps between spaces (including the maps of a covering system) are basepoint-preserving.  Furthermore, all homotopies and homotopy equivalences considered in this paper are assumed to be basepoint-preserving.
\end{remark}

Given a covering system $\X$, we will refer to the corresponding inverse limit $X_\infty$ as a \newword{solenoid}.  For each $\alpha\in\D$, such a solenoid has the structure of a fiber bundle over $X_\alpha$ with totally disconnected fibers.  In particular, $X_\infty$ is equipped with projection maps $p_\alpha \colon X_\infty \to X_\alpha$, which commute with the covering maps $p_{\alpha\beta}$ in the sense that $p_{\alpha\beta}\circ p_\beta = p_\alpha$.

Many of the theorems we use apply only in the case where the inverse limit solenoid $X_\infty$ is compact.  By Tychonoff's Theorem, this occurs if and only if each of the spaces $X_\alpha$ in the covering system is compact.  In this case, we will refer to $\X$ as a \newword{compact covering system}, and the resulting solenoid $X_\infty$ as a \newword{compact solenoid}.

We are especially interested in aspherical solenoids.  Recall that a topological space $X$ is \newword{aspherical} if $\pi_n(X) = 0$ for all $n \geq 2$.  If $X$ is a CW complex, this is equivalent to the requirement that the universal cover of $X$ is contractible. For example, any finite graph is aspherical, as is any closed hyperbolic or Euclidean manifold. The fundamental theorem for aspherical CW~complexes is the following:

\begin{theorem}
Let $X$ and $Y$ be connected CW complexes, with $Y$ aspherical. Then every homomorphism $\pi_1(X) \to \pi_1(Y)$ is induced by a map $X \to Y$ that is unique up to homotopy.
\label{thm:hatcher}
\end{theorem}
\begin{proof}
See \cite[Prop. 1B.9, pg. 90]{hatcher}.
\end{proof}

This theorem has several important consequences in topology, including the Dehn-Nielsen Theorem for closed surfaces.  Our main result is a generalization of Theorem~\ref{thm:hatcher} to compact aspherical solenoids.

The following proposition characterizes aspherical solenoids:

\begin{proposition}
Let\/ $\X$ be a covering system with inverse limit $X_\infty$. Then $X_\infty$ is aspherical if and only if each CW~complex $X_\alpha$ of\/ $\X$ is aspherical.
\label{prop:asphericalsolenoids}
\end{proposition}
\begin{proof}Let $\alpha\in\D$, and recall that the map $p_\alpha\colon X_\infty\to X_\alpha$ is the projection for a fiber bundle with totally disconnected fibers.  Using the long exact sequence of homotopy groups of a fibration, it follows immediately that $\pi_n(X) \cong \pi_n(X_\alpha)$ for all $n\geq 2$.  In particular, $X_\infty$ is aspherical if and only if $X_\alpha$ is aspherical.
\end{proof}

Using this proposition, we can construct several examples of compact aspherical solenoids:

\begin{example}Let $P = \{p_1, p_2, \ldots \}$ be a sequence of prime
numbers, and consider the following covering system:
\[
S^1 \;\overset{p_1}{\longleftarrow}\; S^1
\;\overset{p_2}{\longleftarrow}\; S^1
\;\overset{p_3}{\longleftarrow}\; \cdots\text{.}
\]
Here each space is a circle, and each map is a covering of the indicated degree.  The resulting solenoid is known as the \newword{$\boldsymbol P$-adic solenoid} \cite{gumerov}.  In the special case where $P$ is an infinite sequence of~$2$'s, this is Smale's original \newword{dyadic solenoid}.
\label{ex:padic}
\end{example}

\begin{example}
Let $S$ be a Euclidean or hyperbolic surface, and let $\X$ be a system of finite-sheeted covers of~$S$.  Then the resulting solenoid $X_\infty$ is a \newword{surface solenoid}.  In particular, if $\X$ is a system of tori, then $X_\infty$ is a \newword{torus-like solenoid}~\cite{eda-mandic-matijevic}. If $\X$ is a system of Klein bottles, then $X_\infty$ is a \newword{Klein-bottle-like solenoid}~\cite{matijevic}.

If $\X$ is a system of hyperbolic surfaces, then $X_\infty$ is called a \newword{hyperbolic solenoid}.  In particular, the \newword{unviersal hyperbolic solenoid} is the solenoid obtained when $\X$ is the lattice of all finite-sheeted covers of a hyperbolic surface~$S$.  Because any two hyperbolic surfaces have a common cover (up to homeomorphism), the homeomorphism type of the universal hyperbolic solenoid does not depend on the base hyperbolic surface~$S$.

In~\cite{odden}, Odden proved that the mapping class group of the universal hyperbolic solenoid is precisely the group of virtual automorphisms of the fundamental group of a hyperbolic surface.  This theorem was the starting point for our work, and we recover a version of this result using purely topological methods in Section~8.
\label{ex:surfsol}
\end{example}

\begin{example}Let $\Gamma$ be a finite graph with nonabelian fundamental group, and let $\X$ be a system of finite-sheeted covers of~$\Gamma$.  Then the resulting solenoid $X_\infty$ is a \newword{graph solenoid}.  Such solenoids have not been studied extensively, and we think they merit further attention.

As a special case, if $\X$ is the lattice of all finite-sheeted covers of $\Gamma$, then the resulting solenoid $X_\infty$ is a \newword{universal graph solenoid}.  We will prove in Section~8 that the homotopy type of a universal graph solenoid does not depend on the base graph~$\Gamma$.  Indeed, if we restrict to the case where the base graph~$\Gamma$ is trivalent, then the homeomorphism type of $X_\infty$ is uniquely determined.

As we will show in Section~8, the self-homotopy-equivalences of a universal graph solenoid are in one-to-one correspondence with the virtual automorphisms of a nonabelian free group.  This can be viewed as an analog of Odden's result in the context of automorphisms of free groups and outer space (see~\cite{culler-vogtmann}).
\label{ex:univgraphsol}
\end{example}

\section{Morphisms and limit maps}

It is possible to construct maps between solenoids by lifting maps between spaces of the corresponding covering systems.  This involves the notion of a morphism between covering systems, which we will now define.

Let $\X$ and $\Y$ be covering systems, and let $f\colon X_\mu \to Y_\nu$ be a map. By a \newword{lift} of $f$, we mean a map $X_\alpha \to Y_\beta$ (where $\alpha \geq \mu$ and $\beta \geq \nu$) making the following diagram commute:
    \[
    \xymatrix@=0.3in@M=0.5em{
    X_\alpha\ar@{-->}[r]\ar[d] & Y_\beta\ar[d] \\
    X_\mu\ar[r]_f & Y_\nu
    }
    \]
If a lift of $f$ from $X_\alpha$ to $Y_\beta$ exists, it must be unique \cite[Prop.~1.34]{hatcher}.

A map $f\colon X_\mu \to Y_\nu$ is \newword{liftable} with respect to covering systems $\X$ and~$\Y$ if for every $\beta \geq \nu$, there exists a lift of $f$ with codomain~$Y_\beta$. That is, $f$ is liftable if it can be lifted ``all the way up'' the tower of covers.  A \newword{morphism} from $\X$ to $\Y$ is an equivalence class of liftable maps, where two maps $f\colon X_\mu \to Y_\nu$ and $f'\colon X_{\mu'} \to Y_{\nu'}$ are considered equivalent if they have a common lift. (This is a special case of a more general notion of morphism between inverse systems. See \cite{MaSe}.)

The \newword{composition} of two morphisms $\f\colon \X\to \Y$ and $\g\colon \Y\to \Z$, denoted $\g \circ \f \colon \X\to \Z$, is obtained by composing a representative $g\colon Y_\nu \to Z_\xi$ for $\g$ with a corresponding representative $f\colon X_\mu \to Y_\nu$ for~$\f$.  Under this notion of composition, covering systems and morphisms form a category.

\begin{remark}
It is also possible to define morphisms from a topological space $X$ to a covering system $\Y$.  In this case, a \newword{lift} of a map $f\colon X \to Y_\nu$ is any map $X \to Y_\beta$ (where $\beta \geq \nu$) making the following diagram commute:
    \[
    \xymatrix@C=0.3in@R=0in@M=0.5em{
    & Y_\beta\ar[dd] \\
    X \ar@{-->}[ur]\ar[dr]_f &  \\
    & Y_\nu
    }
    \]
A map is \newword{liftable} if it can be lifted ``all the way up'' the covers of $\Y$, and two maps are equivalent if they have a common lift.  An equivalence class of liftable maps is a \newword{morphism} from $X$ to~$\Y$.

In the case where $X$ is a CW complex, we can view $X$ itself as a covering system $\X$ indexed over a one-element directed set, in which case a morphism $X\to\Y$ is precisely the same thing as a morphism $\X\to\Y$.  However, we are also interested in morphisms whose domains are not CW~complexes.  For example, if $Y_\infty$ is the inverse limit solenoid for~$\Y$, then the projection maps $q_\nu\colon Y_\infty\to Y_\nu$ are all equivalent, and define a morphism $\q\colon Y_\infty\to\Y$.
\end{remark}

We can use the language of morphisms to succinctly state the universal property for the inverse limit of a covering system.  Let $\Y$ be a covering system with inverse limit $Y_\infty$, and let $\q\colon Y_\infty \to \Y$ be the morphism determined by the associated projection maps.  Then for any space $X$ and any morphism $\f\colon X \to \Y$, there exists a unique map $f_\infty\colon X \to Y_\infty$ making the following diagram commute:
    \[
    \xymatrix@=0.4in@M=0.5em{
     & Y_\infty\ar[d]^{\q} \\
    X\ar[r]_{\f}\ar@{-->}[ur]^{f_\infty} & \Y
    }
    \]
The following proposition generalizes this universal property to morphisms between covering systems:

\begin{proposition}
Let\/ $\X$ and $\Y$ be covering systems, and let\/ $\f\colon \X\to\Y$ be a morphism.
Then there exists a unique map $f_\infty\colon X_\infty \to Y_\infty$ between the corresponding solenoids for which the following diagram commutes: \[
    \xymatrix@C=0.35in@R=0.4in@M=0.5em{
    X_\infty\ar@{-->}[r]^{f_\infty}\ar[d]_{\p} & Y_\infty\ar[d]^{\q} \\
    \X\ar[r]_{\f} & \Y
    }
    \]
    \label{prop:lift}
\end{proposition}
\begin{proof}The composition $\f\circ\p \colon X_\infty \to \Y$ is a morphism whose domain is a single space.
Thus, by the universal property of inverse limits, there exists a unique map
$f_\infty\colon X_\infty \to Y_\infty$ such that $\q \circ f_\infty = \f \circ \p$.
\end{proof}

We shall refer to the map $f_\infty$ as the \newword{limit} of the morphism $\f$, and any map that is the limit of some morphism will be called a \newword{limit map}.  Though a general map between two solenoids is not a limit map, we will show in the next section that, under certain conditions, any map between two solenoids is homotopic to a limit map.

\section{Homotopies of limit maps}
In this section, we will state an important theorem from shape theory (obtained from \cite{MaSe}) and examine its consequences for covering systems.  Before proceeding, we should point out that our context differs from the context of \cite{MaSe} in three respects:
\begin{enumerate}
\item In \cite{MaSe}, our definition of a morphism between covering systems is replaced by a more general definition of a morphism between inverse systems.\smallskip
\item In \cite{MaSe}, the ``nice'' spaces used to build inverse systems are any spaces with the homotopy type of an absolute neighborhood retract~(ANR).  This includes the CW~complexes required by our definition of a covering system.\smallskip
\item In \cite{MaSe}, morphisms between inverse systems of topological spaces are only defined up to homotopy.
\end{enumerate}
We would like to expand on point (3).  All of our topological definitions have been using the category Top* of basepointed topological spaces and continuous maps, but \cite{MaSe} uses the category HTop (or HTop*) whose morphisms are homotopy classes of continuous maps.  This leads to a slightly different notion of morphism between two covering systems $\X$ and $\Y$.

In particular, using the category HTop (or HTop*), a map would be considered liftable if it can be lifted up to homotopy, and two liftable maps would be considered equivalent if they have a homotopic pair of lifts.  Fortunately, this makes very little difference in the context of covering systems.  In fact, it follows from the homotopy lifting property of covers \cite[Prop.~1.30]{hatcher} that any map $f \colon X_\mu \to Y_\nu$ that is liftable up to homotopy is in fact liftable.  Therefore, every HTop morphism $\X \to \Y$ is simply a homotopy class of morphisms $\X\to\Y$, where two morphisms are said to be \newword{homotopic} if they have a homotopic pair of representatives.

The following theorem is the primary shape-theoretic result that we will use:

\begin{theorem}Let\/ $\X$ and\/ $\Y$ be covering systems, with\/ $\X$ compact.  Then for any map\/ $F\colon X_\infty \to Y_\infty$, there exists a morphism\/ $\f\colon \X\to \Y$, unique up to homotopy, for which the following diagram commutes up to homotopy:
    \[
    \xymatrix@C=0.25in@R=0.3in@M=0.5em{
    X_\infty\ar[r]^{F}\ar[d] & Y_\infty\ar[d] \\
    \X \ar[r]_{\f} & \Y
    }
    \]
\label{thm:shapetheory}
\end{theorem}
\begin{proof}This is essentially a version of \cite[Thm.~I.5.9, pg.~65]{MaSe}, which states that the morphism $X_\infty \to \X$ is an \newword{HTop-expansion} of $X_\infty$. There are only a few differences between the two theorems:
\begin{enumerate}
\item We are assuming that each space has a basepoint, while \cite{MaSe} does not.  Thus, our theorem technically follows from \cite[Thm.~I.5.13, pg.~71]{MaSe}, which concerns inverse limits in the category of compact pairs.

\smallskip\item The conclusion of the theorem in \cite{MaSe} is that the map $F$ defines a unique HTop morphism $[\f]\colon\X\to \Y$ for which the diagram above commutes up to homotopy.  By the discussion above, this HTop morphism $[\f]$ is simply the homotopy class of some morphism $\f\colon \X\to \Y$.\qedhere
\end{enumerate}
\end{proof}

\noindent We will be using the following two consequences of this theorem:

\begin{theorem}Let\/ $\X$ and\/ $\Y$ be covering systems, with\/ $\X$ compact.
Then any map $X_\infty \to Y_\infty$ is homotopic to a limit map.
\label{cor:fin}
\end{theorem}
\begin{proof}Let $F\colon X_\infty \to Y_\infty$ be a map, and let $\f\colon\X\to\Y$ be the associated morphism given by Theorem \ref{thm:shapetheory}.  If $f\colon X_\alpha \to Y_\beta$ is a representative for~$\f$, then the following square must
commute up to homotopy:
    \[
    \xymatrix@C=0.25in@R=0.3in@M=0.5em{
    X_\infty\ar[r]^{F}\ar[d]_{p_\alpha} & Y_\infty\ar[d]^{q_\beta} \\
    X_\alpha\ar[r]_{f} & Y_\beta
    }
    \]
Now, the composition $f\circ p_\alpha$ is equal to $q_\beta \circ f_\infty$, where $f_\infty \colon X_\infty \to Y_\infty$ is the limit of $\f$. From the above square, we conclude that $q_\beta \circ f_\infty$ and $q_\beta \circ F$ are homotopic.  Since $q_\beta$ is a fibration, it follows that $f_\infty$ is
homotopic to~$F$.
\end{proof}

\begin{theorem}Let\/ $\X$ and\/ $\Y$ be covering systems, with\/ $\X$ compact, and let\/ $\f,\g\colon \X\to\Y$ be morphisms.  Then\/ $\f$ and\/ $\g$ are homotopic if and only if the limit maps $f_\infty, g_\infty\colon X_\infty \to Y_\infty$ are homotopic.
\label{cor:hom}
\end{theorem}
\begin{proof}If $f_\infty$ and $g_\infty$ are homotopic, then $\f$ and $\g$ must also be homotopic by the uniqueness portion
of Theorem \ref{thm:shapetheory}.  For the converse, suppose that $\f$ and $\g$ are homotopic as morphisms,
and let $f,g\colon X_\alpha \to Y_\beta$ be a homotopic pair of representatives for $\f$ and $\g$.  Composing with the projection $p_\alpha\colon X_\infty \to X_\alpha$, we find that the diagonal maps of the following two commutative squares are homotopic:
    \[
    \xymatrix@C=0.25in@R=0.3in@M=0.5em{
    X_\infty\ar[r]^{f_\infty}\ar[d]_{p_\alpha} & Y_\infty\ar[d]^{q_\beta} \\
    X_\alpha\ar[r]_{f} & Y_\beta
    }\qquad\qquad
    \xymatrix@C=0.25in@R=0.3in@M=0.5em{
    X_\infty\ar[r]^{g_\infty}\ar[d]_{p_\alpha} & Y_\infty\ar[d]^{q_\beta} \\
    X_\alpha\ar[r]_{g} & Y_\beta
    }
    \]
Since the projection map $q_\beta \colon Y_\infty \to Y_\beta$ is a fibration, it follows that $f_\infty$ is homotopic to $g_\infty$.
\end{proof}

\section{Pro-groups and filtered groups}
A \newword{pro-group} is an inverse system of groups and homomorphisms, indexed over a directed set.  For example, any inverse sequence
\[
G_0 \;\overset{\rho_1}{\longleftarrow}\; G_1
\;\overset{\rho_2}{\longleftarrow}\; G_2
\;\overset{\rho_3}{\longleftarrow}\; \cdots
\]
is a pro-group, indexed over the natural numbers.

A \newword{filtered group} is a pro-group $\G$ for which each bonding homomorphisms $\rho_{\alpha\beta}\colon G_\beta \to G_\alpha$ is injective.  For example, the lattice of all subgroups of a group $G$ forms a filtered group whose bonding homomorphisms are inclusions.  Many of the standard definitions in the theory of pro-categories are somewhat simpler in the context of filtered objects, and for this reason we will restrict ourselves to filtered groups.

There is a notion of morphism between filtered groups, similar to the morphisms between covering systems developed in Section~3.  If $\G$ and $\H$ are filtered groups and $\varphi\colon G_\mu\to H_\nu$ is a homomorphism, a \newword{restriction} of $\varphi$ is any
homomorphism $G_\alpha \to H_\beta$ (where $\alpha \geq \mu$ and $\beta \geq \nu$) making the following diagram commute:
    \[
    \xymatrix@=0.3in@M=0.5em{
    G_\alpha\ar@{-->}[r]\ar[d] & H_\beta\ar[d]\\
    G_\mu\ar[r]_\varphi & H_\nu
    }
    \]
We say that $\varphi$ is \newword{compatible} with $\G$ and $\H$ if for every $\beta \geq \nu$, there exists a restriction of $\varphi$ with codomain~$H_\beta$.  A \newword{morphism} from $\G$ to $\H$ is an equivalence class of compatible homomorphisms, where two homomorphisms $\varphi\colon G_\mu \to H_\nu$ and $\varphi'\colon G_{\mu'}\to H_{\nu'}$ are considered equivalent if they have a common restriction.  (This is a special case of a more general notion of a morphism between pro-groups.  See \cite{MaSe}.)

The \newword{composition} of morphisms $\Bphi\colon \G\to\H$ and $\Bpsi\colon \H\to\K$, denoted $\Bpsi\circ\Bphi\colon\G\to\K$, is obtained by composing a representative $\psi\colon H_\nu \to K_\xi$ of $\Bpsi$ with a corresponding representative $\varphi\colon G_\mu \to H_\nu$ of $\Bphi$.
Under this notion of composition, filtered groups and morphisms form a category.

A morphism $\Bphi\colon \G \to \H$ is called an \newword{isomorphism} if there exists a morphism $\Bpsi\colon \H\to\G$ for which $\Bpsi \circ \Bphi$ is the identity on $\G$, and $\Bphi \circ \Bpsi$ is the identity on $\H$.  An isomorphism from a filtered group $\G$ to itself is called an
\newword{automorphism}.

\begin{example}
Let $P = (p_1,p_2,\ldots)$ be a sequence of prime numbers, and consider the filtered group
\[
\mathbb{Z} \;\overset{p_1}{\longleftarrow}\; \mathbb{Z} \;\overset{p_2}{\longleftarrow}\; \mathbb{Z} \;\overset{p_3}{\longleftarrow}\;\cdots,
\]
where each homomorphism is multiplication by the indicated prime.  Such a filtered group can be thought of as a descending sequence of subgroups of the integers:
\[
\mathbb{Z} \;\;\geq\;\; m_1\mathbb{Z} \;\;\geq\;\; m_2\mathbb{Z} \;\;\geq\;\; m_3\mathbb{Z} \;\;\geq\;\; \cdots
\]
where $m_k = p_1 \cdots p_k$.  We wish to classify all morphisms between filtered groups of this form.

First, recall that every homomorphism $\varphi\colon m\mathbb{Z}\to n\mathbb{Z}$ between two nontrivial subgroups of the integers has the form
\[
\varphi(k) = \frac{a}{b}\,k
\]
where $a/b \in \mathbb{Q}$ and $(a/b)m$ is an integer multiple of $n$.  Now, suppose that $\G$ and $\H$ are two filtered groups of the above form, say
\[
\begin{array}{rcrcrcrcr}
\G:\quad&\mathbb{Z} &\geq& m_1\mathbb{Z} &\geq& m_2\mathbb{Z} &\geq& m_3\mathbb{Z} &\geq \cdots \\[3pt]
\H:\quad&\mathbb{Z} &\geq& n_1\mathbb{Z} &\geq& n_2\mathbb{Z} &\geq& n_3\mathbb{Z} &\geq \cdots
\end{array}
\]
Then a morphism $\Bphi\colon\G\to\H$ corresponds to a rational number $a/b$ with the following property: for every $j\in\mathbb{N}$, there exists a $i\in\mathbb{N}$ so that $(a/b)m_i$ is an integer multiple of~$n_j$.

We can understand this criterion better using infinite prime factorizations. Let $2^{r_2}3^{r_3}5^{r_5}\cdots$ be the limiting prime factorization of~$m_k$ as~$k\to\infty$, where each exponent $r_k$ is an element of $\mathbb{N}\cup\{\infty\}$.  Similarly, let $2^{s_2}3^{s_3}5^{s_5}\cdots$ be the limiting prime factorization of~$n_k$ as~$k\to\infty$.  Then a rational number $a/b$ corresponds to a morphism $\Bphi\colon\G\to\H$ if and only if
\[
b \;\bigl|\; \bigl(2^{r_2}3^{r_3}5^{r_5}\cdots\bigr) \qquad\text{and}\qquad
\bigl(2^{s_2}3^{s_3}5^{s_5}\cdots\bigr) \;\bigl|\; (a/b)\bigl(2^{r_2}3^{r_3}5^{r_5}\cdots\bigr),
\]
where products and divides should be interpreted in the obvious way for infinite prime factorizations.

Such a morphism is invertible if and only if the reciprocal $b/a$ determines a morphism from $\H$ to $\G$.  Using this criterion, it is easy to verify the following facts:
\begin{itemize}
\item $\G$ and $\H$ are isomorphic if and only if $2^{s_2}3^{s_3}5^{s_5}\cdots$ is a rational multiple of $2^{r_2}3^{r_3}5^{r_5}\cdots$, i.e.~if and only if the sum
    \[
    |r_2 - s_2| + |r_3 - s_3| + |r_5 - s_5| + \cdots
    \]
    is finite.\smallskip
\item A rational number $a/b$ corresponds to an automorphism of $\G$ if and only if $r_p = \infty$ for every prime $p$ dividing either $a$ or $b$.
\end{itemize}
\label{ex:padicalg}
\end{example}

\section{Fundamental pro-groups}

\noindent Armed with the necessary algebraic tools, we are now in a position to describe the algebraic topology of covering systems and solenoids.  If $\X$ is a covering system, the \newword{fundamental pro-group} of $\X$, denoted $\Bpi_1(\X)$, is defined as follows:
\begin{itemize}
\item The groups of $\Bpi_1(\X)$ are the fundamental groups $\pi_1(X_\alpha)$.\smallskip
\item The bonding homomorphisms of $\Bpi_1(\X)$ are the homomorphisms on fundamental group induced by the covering maps of $\X$.
\end{itemize}
Since the maps of $\X$ are coverings, each of the induced homomorphisms is injective, and therefore the fundamental pro-group $\Bpi_1(\X)$ is actually a filtered group.

The following proposition establishes the existence of induced morphisms for fundamental pro-groups of covering systems:

\begin{proposition}
Let\/ $\X$ and\/ $\Y$ be covering systems, let $f\colon X_\mu \to Y_\nu$ be a map, and let $f_*\colon \pi_1(X_\mu)\to\pi_1(Y_\nu)$ be the induced homomorphism.  Then:\smallskip
\begin{enumerate}
\item The map\/ $f$ is liftable with respect to\/ $\X$ and\/ $\Y$ if and only if $f_*$ is compatible with\/ $\Bpi_1(\X)$ and\/ $\Bpi_1(\Y)$. \smallskip
\item If a map $g\colon X_\alpha\to Y_\beta$ is equivalent to $f$, then the induced homomorphisms $f_*$ and $g_*$ are equivalent as well.
\end{enumerate}
\label{prop:liftable}
\end{proposition}
\begin{proof} For (1), if $f$ is liftable with respect to $\X$ and $\Y$, then $f_*$ must be compatible with $\Bpi_1(\X)$ and $\Bpi_1(\Y)$, since the homomorphisms induced by the lifts of $f$ provide the necessary restrictions of~$f_*$.  For the converse, suppose that $f_*$ is compatible, and let $Y_\beta$ be a cover of $Y_\nu$.  Then there must exist a cover $X_\alpha$ of $X_\mu$ so that $f_*$ restricts to a homomorphism $\pi_1(X_\alpha) \to \pi_1(Y_\beta)$.  By the lifting criterion for covers \cite[Prop. 1.33]{hatcher}, it follows that $f$ lifts to a map $X_\alpha\to Y_\beta$, and therefore $f$ is liftable.

Finally, if $g$ is equivalent to $f$, then the homomorphism induced by their common lift is a common restriction of $f_*$ and $g_*$, and therefore $f_*$ and $g_*$ are equivalent as well.
\end{proof}

It follows from this proposition that any morphism $\f\colon\X\to\Y$ gives a well-defined pro-group morphism $\f_*\colon \Bpi_1(\X)\to \Bpi_1(\Y)$, where the representatives for $\f_*$ are the homomorphisms induced by the representatives for~$\f$.  We shall refer to $\f_*$ as the pro-group morphism \newword{induced} by~$\f$.  Using this definition, fundamental pro-group becomes a functor from the category of covering systems to the category of filtered groups.

The next proposition shows that the induced pro-group morphism depends only on the homotopy class of the limit map:

\begin{proposition}Let\/ $\X$ and\/ $\Y$ be covering systems, with\/ $\X$ compact, and let $X_\infty$ and $Y_\infty$ be the inverse limits. Let\/ $\f,\g\colon \X\to\Y$ be morphisms, and let $F,G\colon X_\infty\to Y_\infty$ be the corresponding limit maps.  If $F$ and $G$ are homotopic, then the induced pro-group morphisms\/ $\f_*,\g_*\colon \Bpi_1(\X)\to\Bpi_1(\Y)$ are the same.
\end{proposition}
\begin{proof}
Suppose $F$ and $G$ are homotopic.  By Theorem~\ref{cor:hom}, it follows that the morphisms $\f$ and $\g$ are homotopic.  In particular, we can find a homotopic pair of representatives $f,g\colon X_\mu\to Y_\nu$ for $\f$ and~$\g$.  Then the induced homomorphisms $f_*,g_*\colon \pi_1(X_\mu)\to\pi_1(Y_\nu)$ are equal, and therefore $\f_* = \g_*$.
\end{proof}

Since any map between compact solenoids is homotopic to a limit map (see Theorem~\ref{cor:fin}), this proposition gives us a well-defined homomorphism $F_* \colon \Bpi_1(\X) \to \Bpi_1(\Y)$ induced by \textit{any} map $F\colon X_\infty\to Y_\infty$.

As a consequence, we can now prove that homotopy-equivalent compact solenoids have isomorphic fundamental pro-groups:

\begin{proposition}
Let\/ $\X$ and\/ $\Y$ be compact covering systems, and let $X_\infty$ and $Y_\infty$ be the corresponding solenoids.  Then any homotopy equivalence $F\colon X_\infty \to Y_\infty$ induces an isomorphism $F_*\colon \Bpi_1(\X) \to \Bpi_1(\Y)$.
\label{prop:inducediso}
\end{proposition}
\begin{proof}Let $G\colon Y_\infty \to X_\infty$ be a homotopy inverse for $F$, so $G\circ F$ and
$F\circ G$ are homotopic to identity maps.  Then $G_* \circ F_* = (G\circ F)_*$ is the identity morphism of $\Bpi_1(\X)$ and
$F_* \circ G_* = (F \circ G)_*$ is the identity morphism of $\Bpi_1(\Y)$, so $F_*$ and $G_*$ are an inverse pair of isomorphisms.
\end{proof}

The above proposition tells us that the fundamental pro-group of a compact covering system depends only on the homotopy class of the inverse limit solenoid.  In particular, different compact covering systems with the same inverse limit must have isomorphic fundamental pro-groups.  This lets us make the following definition:

\begin{definition}
If $X$ is a compact solenoid, the \newword{fundamental pro-group} of~$X$, denoted~$\Bpi_1(X)$, refers to the fundamental pro-group of any covering system whose inverse limit is~$X$.
\end{definition}

\begin{remark}
The fundamental pro-group we are using here is the same as the ``fundamental trope'' introduced by Fox in \cite{fox}.  Fox proved a Galois correspondence between overlays of a connected metric space and pro-subgroups of its fundamental pro-group, where an \textit{overlay} is a certain generalization of a covering space.  (This result was recently extended to include all connected topological spaces \cite{mardesic-matijevic}.)  Among other things, these results prove that finite-sheeted covers of a compact solenoid $X$ are in one-to-one correspondence with finite-index pro-subgroups of $\Bpi_1(X)$.
\end{remark}

\section{Compact aspherical solenoids}

We are now ready to prove our main theorem for compact aspherical solenoids:

\begin{theorem}
Let\/ $X$ and\/ $Y$ be compact solenoids, where\/ $Y$ is aspherical. Then every pro-group morphism $\Bpi_1(X) \to \Bpi_1(Y)$ is induced by a map $X \to Y$ that is unique up to homotopy.
\label{thm:nielsen}
\end{theorem}

\begin{proof}[Proof of Theorem \ref{thm:nielsen}]
Let $\X$ and $\Y$ be covering systems whose inverse limits are $X$ and $Y$, respectively.  To prove the existence statement, let $\Bphi\colon \Bpi_1(\X) \to \Bpi_1(\Y)$ be a morphism, and let $\varphi\colon \pi_1(X_\mu) \to \pi_1(Y_\nu)$ be a representative for $\Bphi$.  Since $Y_\nu$ is aspherical, by Theorem~\ref{thm:hatcher} there exists a map $f\colon X_\mu \to Y_\nu$ for which $f_* = \varphi$.  By Proposition~\ref{prop:liftable}, the map $f$ is liftable, so it defines a morphism $\f\colon \X \to \Y$.  The limit $f_\infty \colon X \to Y$ of this morphism exists by Proposition~\ref{prop:lift}.  Then $(f_\infty)_* = \f_* = \Bphi$ which proves the existence portion of the theorem.

For uniqueness, let $F, G \colon X \to Y$ and suppose that $F_* = G_*$.  By Theorem~\ref{cor:fin}, there exist morphisms $\f, \g \colon \X \to \Y$ so that $f_\infty$ is homotopic to $F$ and $g_\infty$ is homotopic to~$G$.  By definition $\f_* = \g_*$, so we can choose representatives $f,g\colon X_\alpha \to Y_\beta$ of $\f$ and $\g$ for which $f_* = g_*$.  Since $Y_\beta$ is aspherical, it follows that $f$ and $g$ are homotopic.  Then $\f$ and $\g$ are homotopic as morphisms, so the limit maps $f_\infty$ and $g_\infty$ must be homotopic by Theorem~\ref{cor:hom}.  We conclude that $F$ is homotopic to $G$, which proves the uniqueness portion.
\end{proof}

We shall consider several applications of Theorem~\ref{thm:nielsen}, beginning with the following:

\begin{theorem} Two compact aspherical solenoids are homotopy equivalent if and only if their fundamental pro-groups are isomorphic.
\label{cor:homotopyclass}
\end{theorem}
\begin{proof}
The forward direction was proven in Proposition~\ref{prop:inducediso}.  For the converse, let $X$ and $Y$ be compact aspherical solenoids, and suppose there exists an isomorphism $\Bphi \colon \Bpi_1(X) \to \Bpi_1(Y)$.  Let  $\Bpsi\colon \Bpi_1(Y) \to \Bpi_1(X)$ be the inverse of $\Bphi$, so both $\Bpsi \circ \Bphi$ and $\Bphi \circ \Bpsi$ are identity morphisms.  Since $X$ and $Y$ are aspherical, it follows from Theorem~\ref{thm:nielsen} that there exist maps $F\colon X\to Y$ and $G\colon Y\to X$ for which $F_* = \Bphi$ and $G_* = \Bpsi$.  Then $(F\circ G)_*$ and $(G\circ F)_*$ are both identity morphisms, and hence by Theorem~\ref{thm:nielsen} both $F \circ G$ and $G\circ F$ are homotopic to the identity.  Thus $F$ is the desired homotopy equivalence.
\end{proof}

\begin{example}
Theorem \ref{thm:nielsen} can be used to classify all maps between $P$-adic solenoids up to homotopy.  Given a sequence $P = \{p_1,p_2,\ldots\}$
of prime numbers, the fundamental pro-group of the corresponding $P$-adic solenoid has the form
\[
\mathbb{Z} \;\overset{p_1}{\longleftarrow}\; \mathbb{Z} \;\overset{p_2}{\longleftarrow}\; \mathbb{Z} \;\overset{p_3}{\longleftarrow}\;\cdots,
\]
where each homomorphism is multiplication by the indicated prime.  Morphisms between filtered groups of this form were discussed in Example~\ref{ex:padicalg}, and by Theorem~\ref{thm:nielsen} each such morphism is induced by a unique homotopy class of maps.

Similarly, Theorem~\ref{cor:homotopyclass} can be applied to $P$-adic solenoids, as it implies that two $P$-adic solenoids are homotopy equivalent if and only if they have isomorphic fundamental pro-groups.  Isomorphisms between filtered groups of this form were also discussed in Example~\ref{ex:padicalg}, producing a classification of $P$-adic solenoids up to homotopy equivalence.

These classifications are well-known in the case of $P$-adic solenoids.  See \cite{kwapisz} for a detailed discussion.
\label{ex:padicclass}
\end{example}

Given a space $X$, the \newword{homotopy self-equivalence group} of $X$, denoted $\mathcal{E}(X)$, is the group of all homotopy classes of (basepoint-preserving) homotopy equivalences $X \to X$.  This group has been studied extensively for various spaces \cite{arkowitz} \cite{kahn}.  The following theorem characterizes this group for compact aspherical solenoids.

\begin{theorem} Let $X$ be a compact aspherical solenoid.  Then the group $\mathcal{E}(X)$ of homotopy self-equivalences of~$X$ is isomorphic to the group of automorphisms of the pro-group\/~$\Bpi_1(X)$.
\label{cor:mappingclassgroup}
\end{theorem}
\begin{proof}
By Theorem~\ref{thm:nielsen}, the semigroup of all homotopy classes of maps $X\to X$ is isomorphic to the semigroup of
all morphisms $\Bpi_1(X) \to \Bpi_1(X)$.  The groups $\mathcal{E}(X)$ and $\mathrm{Aut}\bigl(\Bpi_1(X)\bigr)$ constitute the invertible elements of these semigroups.
\end{proof}

\begin{example}
Let $X$ be a $P$-adic solenoid. By Theorem~\ref{cor:mappingclassgroup}, the homotopy self-equivalence group of $X$ is isomorphic to the automorphism group of the filtered group $\Bpi_1(X)$.  In this case, it is not hard to determine the isomorphism type of $\mathrm{Aut}\bigl(\Bpi_1(X)\bigr)$.  In particular, if $\mathcal{P}$ is the set of primes that appear infinitely often in the sequence $P$, then $\mathrm{Aut}\bigl(\Bpi_1(X)\bigr)$ is isomorphic to the multiplicative group of the nonzero rational numbers whose denominators involve only primes from $\mathcal{P}$.  (This characterization of $\mathcal{E}(X)$ is well-known.  See \cite{kwapisz}.)

For example, if $P$ is an infinite sequence of $2$'s (so that $X$ is the standard dyadic solenoid), then $\mathcal{E}(X)$ is isomorphic to the multiplicative group of the nonzero dyadic rationals. If instead each prime number appears infinitely often in $P$, then $\mathcal{E}(X)$ is isomorphic to the entire multiplicative group of the nonzero rationals.
\end{example}

\section{Virtual automorphisms}

In this section we discuss some applications of our main results to virtual automorphisms of groups.  Recall the following definition:

\begin{definition}
Let $G$ be a group.
\begin{enumerate}
\item A \newword{virtual automorphism} of $G$ is an isomorphism $\varphi\colon H_1\to H_2$, where $H_1$ and $H_2$ are finite-index subgroups of~$G$.\smallskip
\item Two virtual automorphisms are \newword{equivalent} if they agree on some finite-index subgroup of~$G$.\smallskip
\item If $\varphi\colon H_1\to H_2$ and $\psi\colon K_1\to K_2$ are virtual automorphisms of $G$, the \newword{composition} of $\varphi$ and $\psi$ is the virtual automorphism $\psi\circ\varphi$ with domain $\varphi^{-1}(K_1\cap H_2)$ and range $\psi(K_1\cap H_2)$.\smallskip
\item The set of all equivalence classes of virtual automorphisms of $G$ forms a group under composition.  This is the \newword{virtual automorphism group} (or \newword{abstract commensurator}) of $G$, and is denoted~$\VAut(G)$.
\end{enumerate}
\end{definition}

See \cite{farb-handel}, \cite{ivanov}, and \cite[Chapter~6]{zimmer} for more discussion of virtual automorphism groups.  Virtual automorphisms can be studied in the context of pro-groups, using the following definition:

\begin{definition}Let $G$ be a group.  The \newword{virtual core} of $G$, denoted $\mathbf{Virt}(G)$, is the filtered group consisting of all finite-index subgroups of~$G$, with inclusions as bonding homomorphisms.
\end{definition}

\begin{proposition}Let $G$ and $G'$ be groups.
\begin{enumerate}
\item If $H$ is a finite-index subgroup of $G$, then any homomorphism $\varphi\colon H\to G'$ is compatible with $\mathbf{Virt}(G)$ and $\mathbf{Virt}(G')$.
\smallskip

\item A pro-group morphism $\Bphi\colon\mathbf{Virt}(G)\to\mathbf{Virt}(G')$ is an isomorphism if and only if it can be represented by an isomorphism $\varphi\colon H\to H'$, where $H$ is some finite-index subgroup of~$G$ and $H'$ is some finite-index subgroup of~$G'$.
\end{enumerate}
\label{prop:vaut}
\end{proposition}

\begin{remark}
Neither direction of part (2) holds in general for isomorphisms between filtered groups.
\end{remark}

\begin{proof}For statement~(1), let $\varphi\colon H\to G'$ be a homomorphism, and let $H'$ be a finite-index subgroup of~$G'$.  Then $\varphi^{-1}(H')$ is a finite-index subgroup of $H$, and is hence a finite-index subgroup of~$G$. Thus $\varphi$ restricts to a homomorphism from $\varphi^{-1}(H')$ to $H'$.

For statement (2), let $\Bphi\colon\mathbf{Virt}(G) \to \mathbf{Virt}(G')$ be an isomorphism, and let $\Bpsi\colon \mathbf{Virt}(G')\to\mathbf{Virt}(G)$ be the inverse of~$\Bphi$.  Let $\varphi\colon H \to G'$ be a representative for~$\Bphi$, where $H$ is a finite-index subgroup of $G$, and let $\psi\colon H'\to G$ be a representative for~$\Bpsi$, where $H'$ is a finite-index subgroup of~$G'$.

First, consider the composition $\varphi\circ\psi\colon \psi^{-1}(H) \to G'$.  Since $\Bphi\circ\Bpsi$ is the identity, this homomorphism must restrict to the identity map on some finite-index subgroup $K'$ of $G'$.  In particular, $K'$ is contained in the image of~$\varphi$, and therefore the image of $\varphi$ has finite index in~$G'$.

Next, consider the composition $\psi\circ\varphi\colon \varphi^{-1}(H') \to G$.  Since $\Bpsi\circ\Bphi$ is the identity, this homomorphism must restrict to the identity map on some finite-index subgroup $K$ of $G$.  It follows that $\varphi$ is one-to-one on~$K$.  Furthermore, since $\varphi(H)$ has finite index in $G'$, and $K$ has finite index in~$H$, the image $\varphi(K)$ of $K$ must also have finite index in $G'$.  Then $\varphi\colon K \to \varphi(K)$ is an isomorphism that represents~$\Bphi$.

For the converse, suppose that $\Bphi\colon \mathbf{Virt}(G) \to \mathbf{Virt}(G')$ is a pro-group morphism that can be represented by an isomorphism $\varphi\colon H \to H'$, where $H$ is a finite-index subgroup of~$G$ and $H'$ is a finite-index subgroup of~$G'$.  Then the inverse $\varphi^{-1}\colon H'\to H$ is compatible by statement~(1), and therefore defines an inverse morphism $\Bphi^{-1}\colon \mathbf{Virt}(G') \to \mathbf{Virt}(G)$.
\end{proof}

One consequence of this proposition is that the virtual cores $\mathbf{Virt}(G)$ and $\mathbf{Virt}(G')$ of two groups are isomorphic if and only if $G$ and $G'$ are virtually isomorphic, i.e.~if and only if $G$ and $G'$ have a pair of isomorphic subgroups of finite index.

The relationship between virtual automorphisms and virtual cores is given by the following corollary:

\begin{corollary}Let $G$ be a group, and let\/~$\mathbf{Virt}(G)$ be the virtual core of~$G$.  Then the automorphism group of\/~$\mathbf{Virt}(G)$ is isomorphic to the group\/~$\VAut(G)$ of virtual automorphisms of~$G$.\hfill\qedsymbol
\end{corollary}

We now move to topology.  The following definition gives us a universal solenoid over any compact CW~complex:

\begin{definition}
If $X$ is a compact CW complex, the \newword{universal compact solenoid} over $X$, denoted $\mathrm{USol}(X)$ is the solenoid obtained as the inverse limit of the lattice of all finite-sheeted covers of $X$.
\end{definition}

For example, the universal compact solenoid over a hyperbolic surface is the universal hyperbolic solenoid (see Example~\ref{ex:surfsol}).  Similarly,  the universal compact solenoid over a finite graph with nonabelian fundamental group is a universal graph solenoid (see Example~\ref{ex:univgraphsol}).

By the Galois correspondence for covers \cite[Thm.~1.38]{hatcher}, the fundamental pro-group of $\mathrm{USol}(X)$ is precisely the virtual core of~$\pi_1(X)$.  Combining this with Theorem~\ref{thm:nielsen} and Proposition~\ref{prop:vaut}, we obtain:

\begin{theorem}Let $X$ be a compact aspherical CW~complex.  Then the group of homotopy self-equivalences of\/ $\mathrm{USol}(X)$ is isomorphic to the group of virtual automorphisms of $\pi_1(X)$.\hfill\qedsymbol
\label{thm:universalsolenoid}
\end{theorem}

\begin{remark}
This is related to a result of Odden \cite{odden}, who proved a similar statement in the case where $X$ is a hyperbolic surface.  Specifically, Odden proved that the group of isotopy classes of self-homeomorphisms of the universal hyperbolic solenoid is isomorphic to the virtual automorphism group of $\pi_1(X)$.

Though our methods primarily give information about the group of homotopy self-equivalences, in the special case of the universal hyperbolic solenoid they can be used to show that this group is equal to the group of homotopy classes of self-homeomorphisms.  In particular, any homotopy self-equivalence of $\mathrm{USol}(X)$ can be represented by a homotopy equivalence between two finite-sheeted covers of~$X$.  By the Dehn-Nielsen Theorem, such a homotopy equivalence is homotopic to a homeomorphism, and this homeomorphism can then be lifted to yield a self-homeomorphism of~$\mathrm{USol}(X)$.

This observation, together with Theorem~\ref{thm:universalsolenoid}, recovers a version of Odden's result, except that we are using homotopy classes instead of isotopy classes.  It is not obvious how to deduce information about isotopy using our methods.
\end{remark}

Finally, we wish to use our machinery to obtain some information about universal graph solenoids:

\begin{proposition} Let\/ $\Gamma$ and\/ $\Gamma'$ be finite graphs with nonabelian fundamental groups.  Then the universal graph solenoids\/ $\mathrm{USol}(\Gamma)$ and\/ $\mathrm{USol}(\Gamma')$ are homotopy equivalent.
\end{proposition}
\begin{proof}
The fundamental groups $\pi_1(\Gamma)$ and $\pi_1(\Gamma')$ are free groups of finite rank.  By the Schreier Index Formula, any two finite-rank nonabelian free groups are virtually isomorphic.  In particular, the pro-groups $\Bpi_1\bigl(\mathrm{USol}(\Gamma)\bigr)\cong \mathbf{Virt}\bigl(\pi_1(\Gamma)\bigr)$ and $\Bpi_1\bigl(\mathrm{USol}(\Gamma')\bigr)\cong \mathbf{Virt}\bigl(\pi_1(\Gamma')\bigr)$ are isomorphic.  Therefore, by Theorem~\ref{cor:homotopyclass}, the solenoids $\mathrm{USol}(\Gamma)$ and $\mathrm{USol}(\Gamma')$ are homotopy equivalent.
\end{proof}

Indeed, if we restrict to trivalent graphs, then the homeomorphism type of $\mathrm{USol}(\Gamma)$ is uniquely determined.  This follows from the fact that any two trivalent graphs have a common cover.

For any universal graph solenoid $\mathrm{USol}(\Gamma)$, Theorem~\ref{thm:universalsolenoid} says that the group of homotopy self-equivalences of $\mathrm{USol}(\Gamma)$ is isomorphic to the virtual automorphism group of a nonabelian free group.  See \cite{bartholdi-bogopolski} for more information about this group.

\bibliographystyle{plain}

\end{document}